\documentclass[12pt]{article}
\usepackage{amsmath}
\usepackage{amsfonts}
\usepackage{amssymb}
\usepackage{a4wide}

\def\qed{\hfill \square \ }
\def\Frob{\mathrm{Frob}}
\def\p{\mathfrak{p}}
\def\ba{\overline{a}}

\renewcommand\mod{\ \mathrm{mod} \ }

\def\GL{\mathrm{GL}}
\def\Qbar{\overline{\Q}}
\def\Gal{\mathrm{Gal}}
\def\F{\mathbf F}
\def\Fbar{\overline{\F}}

\def\rhobar{\overline{\rho}}
\def\Or{\mathcal O}

\def\Q{\mathbf Q}
\def\Z{\mathbf Z}

\def\GL{\mathrm{GL}}

\newtheorem{theorem}{Theorem}[section]
\newtheorem{lemma}[theorem]{Lemma}

\newtheorem{conj}{Conjecture}

\begin{document}

\author{Frank Calegari\footnote{Supported in part by
the American Institute of Mathematics}}
\title{Mod $p$ representations on Elliptic Curves}
\maketitle

\abstract{Modular Galois representations
$\rhobar: \Gal(\Qbar/\Q) \rightarrow \GL_2(\F_p)$ with cyclotomic
determinant arise from elliptic curves for small $p$.
We show that $\rhobar$ does not necessarily arise from an
elliptic curve whose conductor is as small as possible outside
$p$.
For $p=3$ this disproves a conjecture
of J.~Lario and A.~Rio \cite{Rio}.}
\section{Introduction}

Let $E/\Q$ be an elliptic curve. For any prime number $p$, the
$p$-torsion $E[p]$ is a Galois module that gives rise to a continuous
Galois
representation:
$$\rhobar: \Gal(\Qbar/\Q) \rightarrow \GL_2(\F_p).$$
It follows from standard properties of elliptic curves
(see \cite{AOEC}) that $\rhobar$ is unramified outside $p$ and primes
dividing the conductor $N_E$ of $E$, and that the composition of
$\rhobar$ with the determinant map to $\F^{\times}_p$ is the mod $p$
reduction of the cyclotomic character.  
Conversely, one expects (\cite{Serre}, at least if 
$\rhobar$ is irreducible) that such a $\rhobar$ arises in the usual
way from a modular form $f$ of weight $N(\rhobar)$ and
level $k(\rhobar)$, for certain prescribed $N$ and $k$
(referred to as the Serre level and weight, respectively).
It is not true that $\rhobar$ must arise from an elliptic
curve, however, unless $p$ is small.

\begin{theorem} \label{theorem:genuszero}
Let $p \in \{2,3,5\}$. If $\rhobar: \Gal(\Qbar/\Q) 
\rightarrow \GL_2(\F_p)$ is a modular representation with cyclotomic
determinant, then $\rhobar$ arises from the $p$-torsion of an elliptic
curve.
\end{theorem}

A succinct proof of this result is provided in \cite{Rubin}.
The result follows (not entirely formally) from the fact that
$X(p)$ has genus zero for such $p$.  In this paper, we address
the question of whether the elliptic curve $E$ whose
existence is guaranteed by Theorem~\ref{theorem:genuszero}
can be chosen to have ``minimal'' conductor 
(for a more precise statement, see Theorem~\ref{theorem:small}).
A conjecture along these lines for $p=3$ is made
in \cite{Rio}, and one of the main motivations for this
paper is to find a counterexample to this 
conjecture.
As an afterthought, we discuss some issues related to
representations $\rhobar$ with 
$p \ge 7$. 

\medskip

I would like to thank Akshay Venkatesh and Kiran Kedlaya for
an informative discussion about cubic fields, Ken Ribet for
bringing to my attention the paper~\cite{Rio}, and William
Stein for some computational help.
I would
especially like to thank John Cremona, without whom the main computation
of this paper would not have been possible.

\section{Small $p$}

Let $p \in \{2,3,5\}$, and let $\rhobar: G_{\Q} \rightarrow \GL_2(\F_p)$
be an absolutely irreducible Galois representation arising from an elliptic
curve $E$. If $N(\rhobar)$ denotes the Serre conductor of $E$,
then \emph{a priori} one knows that $N(\rhobar)$ divides $N_E$, the conductor
of $E$. By definition, however, the Serre conductor is coprime to $p$.
Thus if
$\rhobar$ is not finite flat at $p$ (and so $E$ has bad reduction at $p$)
we cannot hope to have an equality $N(\rhobar) = N_E$. Allowing
for this possibility, we may ask 
(given $\rhobar$) whether there exists an
elliptic curve $E$ giving rise to $\rhobar$ such that
$N_E = p^n N(\rhobar)$ for some $n$.
Our main result  is as follows:

\begin{theorem} Let $p \in \{2,3,5\}$. There exists a surjective modular
representation:
$$\rhobar: \Gal(\Qbar/\Q) \rightarrow \GL_2(\F_p)$$
with determinant equal to the cyclotomic character such that
$\rhobar$ does not arise from any elliptic curve of conductor
$p^n N(\rhobar)$ for any $n$, where $N(\rhobar)$ is the Serre
level of $\rhobar$.
\label{theorem:small}
\end{theorem}

When $p=3$, the example we construct 
 provides a counterexample to the following conjecture~\cite{Rio}:

\begin{conj}[Lario, Rio] Let $\mathbf{P}\rhobar: \Gal(\Qbar/\Q) \rightarrow \mathrm{PGL}_2(\F_3)$
be an irreducible representation. Assume that
$\mathbf{P} \rhobar$ has a linear lifting $\rhobar$ to $\GL_2(\F_3)$ with cyclotomic determinant.
Then there is a linear lifting $\rho_{E,3}$ where $E/\Q$ is an elliptic curve having
conductor a power of $3$ times $N(\mathbf{P} \rhobar)$, where $N(\mathbf{P} \rhobar)$
is the minimal Serre conductor of all such liftings.
\end{conj}

Let $\rhobar$ be the representation constructed
for $p=3$ in the proof of theorem~\ref{theorem:small}. 
Then  $N(\rhobar) = 353$ is
prime,
and thus
$N(\mathbf{P} \rhobar) = 353$. Different linear liftings
of $\mathbf{P} \rhobar$ with cyclotomic determinant 
differ by a character $\chi$ with
$\chi^2 = 1$, or equivalently by a quadratic character. 
The conjecture guarantees
the existence of an elliptic curve $E/\Q$ with $\rho_{E,3} = 
\rhobar \otimes \chi$,
and conductor a power of $3$ times $353$. 
In particular, the character $\chi$
can only be ramified at three (if it was ramified at $353$, the Serre conductor
of $\rhobar$ would be divisible by $353^2$). 
If we let $E'$ denote the quadratic twist of $E$ by
$\chi$ then $E'$ will also therefore  have conductor $353$ times a
power of $3$.
On the other hand, by construction, $\rho_{E',3} = \rhobar$,
contradicting the fact that $\rhobar$ does not
arise from an elliptic curve of such a conductor.
Thus the conjecture is false.

\subsection{$p = 2$}

Given a Galois representation $\rhobar$ with image
 $\GL_2(\F_2) \simeq S_3$, there is an obvious way to
construct an associated elliptic curve: any $S_3$ field $L/\Q$ is
the splitting field of an irreducible cubic polynomial $g(x)$, and  the 
elliptic curve $y^2 = g(x)$ gives rise to $\rhobar$.
Let $K$ be a cubic field inside $L$, and let $F$ be the unique quadratic
subfield of $L$.
The Serre weight and level can easily be computed from the
arithmetic of $K$ (of course, $L$ is determined from $K$). In
particular, an odd prime $p$ divides $N(\rhobar)$ if and only if
it divides the field discriminant $\Delta_K$.
Let $V$ be the 
Galois module corresponding to the representation
$\rhobar$.
Let $\alpha \in \Or_K$, and let $f(x)$ be the minimal 
polynomial of $\alpha$. Then $y^2 = f(x)$ is an elliptic
curve whose Galois module $E[2]$ is isomorphic to $V$.
 Any such
elliptic curve arises  from such an $\alpha$.
Moreover, if $\alpha$ is not equal to $c + d^2 \beta$
for some $c$, $d \in \Z$ and $\beta \in \Or_K$ 
the equation $y^2 = f(x)$ provides a
minimal model for
$E$ over $\Z[\frac{1}{2}]$. In particular, 
the  ramification of $E$ at odd primes $\ell$ can
be determined directly from properties of $f(x)$.
When does an $\alpha$
give rise to an elliptic curve with Serre conductor
$2^m \cdot N(\rhobar)$? The polynomial discriminant of $f(x)$
is equal to $\Delta_K$ times the square of the index of
$\Z[\alpha]$ inside $\Or_K$. Moreover the minimal discriminant
of the elliptic curve $y^2 = f(x)$ is equal (up to a power
of two) to the polynomial
discriminant of $f(x)$.
If the prime to $2$ part of $N_E$ is equal to $N(\rhobar)$,
then $E$ has good reduction at every odd prime not dividing
$N(\rhobar)$. 
Thus necessarily  
the polynomial discriminant  of $f(x)$ is not divisible by any
primes other than those already dividing $2 \Delta_K$.  This
is not sufficient, however, since (for example) the cubic
$x^3 - 26$ has discriminant $-2^2 \cdot 3^3 \cdot 13^2$,
and yet the elliptic curve $y^2 = x^3 - 26$ has conductor
$2^6 \cdot 3^2 \cdot 13^2 = 2^6 \cdot 3 \cdot N(\rhobar)$.
If the index of $\Z[\alpha]$ inside $\Or_K$ 
is an exact power of two, however, then the odd part of
the conductor of $E$ is equal to $N(\rhobar)$. 
We prove the following:

\begin{theorem} Let $K/\Q$ be the cubic field
determined by the polynomial $u^3 - u^2 - 2u + 27 = 0$.
Then $\Delta_K = -2063$. Let $L$ be the Galois closure of
$K$, and $\rhobar$ the $\GL_2(\F_2)$ representation
that factors through $\Gal(L/\Q)$. Then $\rhobar$ does not
arise from an elliptic curve of conductor $2^m \cdot 2063$.
Moreover, $K$ is the smallest cubic field
$($with respect to discriminant$)$ with this property.
\end{theorem}

Let $N = N(\rhobar)$. To show that $\rhobar$ does not arise
from an elliptic curve, it suffices to prove that 
there does not exist an element $\alpha \in \Or_K$ such that
$$[\Z \alpha: \Or_K] \in \Z[1/2N]^{\times}.$$
First, however, we eliminate all cubic fields with smaller
discriminant. As we have noted, for such fields it suffices
to construct an element $\alpha \in \Or_K$ whose index is
a power of two.
From the Bordeaux Tables~\cite{Bordeaux}, one can determine
all cubic fields $K/\Q$ with $|\Delta_{K}| \le 2063$. The only
such fields listed whose generating element does not already have
index $1$ or $2$ correspond to the following discriminants
(all from complex cubic fields): $\Delta_K = -1356$, $-1599$, $-1691$, $-1751$,
$-1967$, $-2028$.
These fields do in fact
have elements of index $2$, $1$, $2$, $8$, $1$ and $2$ respectively.
Note that for $\Delta_K = -1751$, there is an element of index $17$ which
also corresponds to an elliptic curve with conductor $2^m \cdot 1751$. The
following table gives these examples, where as usual, $[a_1,a_2,a_3,a_4,a_6]$
denote the elliptic curve
$$y^2 + a_1 y x + a_3 y = x^3 + a_2 x^2 + a_4 x + a_6.$$
\begin{center}
\begin{tabular}{|c|c|c|c|c|}
\hline
$\Delta_K$ & $E$ & $[\Z \alpha: \Or_K]$ & $N(\rhobar)$ & $N_E$  \\
\hline
$-1356$ & $[0,1,0,-9,-21]$ & $2$ & $3 \cdot 113$  & $2^3 N(\rhobar)$ \\
\hline
$-1599$ & $[0,1,0,-14,-27]$ & $1$ & $3 \cdot 13 \cdot 41$ & $2^2 N(\rhobar)$ \\
\hline
$-1691$ & $[0,0,0,-13,-24]$ & $2$ & $19 \cdot 89$ & $2^5 N(\rhobar)$ \\
\hline
$-1751$ & $[0,-6,0,-136,-408]$ & $8$ & $17 \cdot 103$ & $2^6 N(\rhobar)$ \\
        & $[0,0,0,29,-123]$ & $17$ &  $17 \cdot 103$ & $2^5 N(\rhobar)$ \\
\hline
$-1967$ & $[0,-3,0,-16,51]$ & $1$ & $1967$ & $2^4 N(\rhobar)$ \\
\hline
$-2028$ & $[0,-1,0,-17,-27]$ & $2$ & $3 \cdot 13^2$ & $2^3 N(\rhobar)$ \\
\hline
\end{tabular}
\end{center}
 
Thus it suffices to consider  the cubic field of (prime) discriminant
$\Delta_K = -2063$.
Let $K = \Q(u)$, where $u$ satisfies the equation
$$u^3 - u^2 - 2u + 27 = 0.$$
A calculation with {\tt pari} shows that
$$\Or_K = \Z \oplus \left(\Z \cdot u \right) \oplus 
\left( \Z \cdot \frac{u^2 + u}{3} \right).$$
If $\theta = x u + y (u^2 + u)/3$, then the index 
$[\Z \theta: \Or_K]$ is given by the absolute value of the index form, which in this
case is equal to
$$f(x,y) = 3 x^3 + 5 x^2 y + 2 x y^2 + 3 y^3.$$
It now suffices to prove that the equation
$f(x,y) = \pm \  2^m 2063^n$ has no integral solutions.
Without loss of generality we may assume that $x$ and $y$ are co-prime.
Suppose that $m > 0$. Then $f(x,y)$ is even. A simple congruence
check implies that $f(x,y)$ is odd whenever at least one of $x$
or $y$ is odd. Thus $x$ and $y$ are both even, which is impossible
if they are coprime, and hence $m = 0$. Further, for all $x$ and $y$,
$f(x,y) \equiv 0,3,4,5,6 \mod 9$ whereas
 $\pm 2063^n \equiv 1,2,7,8 \mod 9$ for  $3 \nmid n$.
Thus $3$ divides $n$.
 We are therefore reduced to finding elements
of $\Or_K$ of index exactly $2063^{3n}$. Given such an element $\alpha$
its minimal polynomial will have discriminant exactly
$-2063 (2063)^{6n}$. After subtracting perhaps some multiple
of $1/3$ from $\alpha$ (which does not affect the discriminant) the
minimal polynomial of $\alpha$ is
$x^3 - 27 c_4 x - 54 c_6$, where $c_4$, $c_6 \in \Z[1/6]$.
Evaluating the discriminant we find that
$$2^2 3^{9} (c^3_4 - c^2_6) = -  2063 \cdot  (2063)^{6n}.$$
Thus $[324 c_4/2063^{2n},5832 c_6/2063^{3n}]$  
is a $\Z[1/(2 \cdot 3 \cdot 2063)]$
 integral point on the elliptic curve
$$Y^2 = X^3 + 2^4 \cdot 3^3 \cdot 2063 = X^3 + 891216.$$
Using Cremona's program
{\tt mwrank} \cite{Cremona2}, we compute that this curve has no rational points
other than $\infty$, and thus we are done.

\medskip Note that $\rhobar$ does of course arise from the $2$-torsion
of \emph{some} elliptic curve, and moreover $E$ can be chosen such that any
prime $\ell \ne 2,2063$ does not divide the conductor of $E$, as the examples:
$E = [0,0,0,-43,-117]$, $F = [0,-1,0,-2,27]$ of conductors
$2^3 \cdot 5 \cdot 2063$ and $2^4 \cdot 3 \cdot 2063$ demonstrate.
Note also  that we needed to go to a cubic field of rather large discriminant
before we found an $S_3$-representation that did not come from an
elliptic curve of minimal level. We feel this is explained by the
``law of small numbers''. In particular, cubic fields of small discriminant
tend to be quite special, and tend to have integral elements of
very small index. This is not a pattern that persists, however, and
one would expect the example constructed above is the norm rather
than the exception. It is also the reason why we suspected the
conjecture in~\cite{Rio} was false, and set about finding a counterexample.

\subsection{$p = 3$}

This is the case that requires the most computational power,
and I am indebted to John Cremona for reinstalling
and reconfiguring his programs on a 64-bit machine provided
to Harvard by {\tt Sun Microsystems}. 
For reasons analogous to the situation for $p=2$, we may expect
that mod $3$ representations of small Serre conductor do arise from
elliptic curves of small conductor. Unfortunately, one does not have 
fine control
over the set of elliptic curves with fixed mod $3$ representation
in quite the same way as one does for mod $2$ representations.
Thus in
 order to find a candidate mod $3$ representation
that does not come from an elliptic curve, we follow the following
algorithm:
\begin{enumerate}
\item Using William Stein's tables \cite{William}, find all
modular representations of weight $2$ and level $N$ and $3N$. By Serre's
conjecture, any irreducible mod $3$ representation with cyclotomic determinant
and Serre conductor $N$ should arise at these levels.
\item Using Cremona's tables, \cite{Cremona1}, \cite{Cremona}, determine if these
representations come from an elliptic curve of conductor $3^k N$,
for $3^k N \le 20000$ ($20000$ is the current limit of these tables).
If $3^5  N < 20000$ and there are no such elliptic curves
then one is done,
since
$3^5$ is the largest possible power of $3$ dividing the conductor of
an elliptic curve.
\item If the candidate $N$ is larger than
$20000/3^5 \simeq 82.3$ and there
are no elliptic curves of conductor $3^k N < 20000$
giving rise to $\rhobar$, try and construct elliptic
curves of conductor $3^k N$ with large $k$ by 
computing
$\Z[1/6N]$-integral points on the curves
$y^2 = x^3 - 3^k \prod_{\ell | N} \ell^{k_i}$.
This method sometimes enables one to eliminate
$\rhobar$ without having to compute all the
elliptic curves of conductor $3^k N$.
\item Once a representation $\rhobar$ is found that is not
eliminated by any of the previous steps, run Cremona's modular
symbols algorithm \cite{Cremona1} for \emph{all} $3^k N$ 
for $k \le 5$, and determine whether or not those
elliptic curves give rise to $\rhobar$.
\end{enumerate}

To simplify the step $3$, one may choose $N$
to be prime, which significantly cuts down the
number of Diophantine equations one needs to consider.
Note that steps two and three are ultimately
not required for the proof, but are rather present to narrow
down potential examples, since step 4 is very computationally
intensive.
Using this method we find:

\begin{theorem} Let $E$ be the elliptic curve
$[1, 1, 0, -22, -812]$ of conductor $2 \cdot 3 \cdot 353$.
Let $\rhobar: \Gal(\Qbar/\Q) \rightarrow \GL_2(\F_3)$ be the
representation induced from the $3$-torsion of $E$.
The representation $\rhobar$ is surjective, and
the Serre conductor $N(\rhobar) = 353$.
Then $\rhobar$ does not arise from
any elliptic curve of conductor $3^k \cdot 353$.
\end{theorem}

Since $E$ is semistable, the mod $3$ representation is
either reducible or surjective, and a quick check eliminates
the first possibility. Since $E$ is semistable at $2$, we
may check that $E[3]$ is unramified at $2$ by considering
the $2$-adic valuation of the minimal discriminant.
The minimal discriminant is $\Delta_E = -2^{18} \cdot 3 \cdot 353$,
and since $3 | 18$ we conclude that $N(\rhobar) = 353$.
Specifically it arises from a modular form of weight
$2$ and level $3 \cdot 353 = 1059$ (in this case
coming from part of the $3$ torsion on
a modular abelian variety $A_f$ of dimension $17$).

\medskip

\begin{Proof}
It suffices to find all elliptic curves of
conductor $3^k \cdot 353$ for $k = 0, \ldots 5$
and show that none of them give rise to the
mod $3$ representation associated to $E[3]$. This
follows from the next two tables, which give the
trace of Frobenius under the image of $\rhobar$, and
the first few $\ba_p = a_p \mod 3$ for the 
relevant elliptic curves. 

\begin{center}
\begin{tabular}{|c|c|c|c|c|c|c|c|c|c|c|c|c|}
\hline
 $p$ & $2$ & $3$ & $5$ & $7$ & $11$& $13$ & $17$ & $19$ 
& $23$ & $29$ & $31$ & $37$ \\ 
\hline
$\mathrm{Trace}(\rhobar(\mathrm{Frob}_p))$ & $0$ &  &
$1$ & $2$ & $1$ & $1$ & $0$ & $1$ & $0$ & $0$ & $1$ & $0$ \\
\hline
\end{tabular}
\end{center}

\begin{center}
\begin{tabular}{|c|l|c|c|c|c|c|c|c|c|c|}
\hline
$N$ & $E$ & $\ba_2$ & $\ba_3$ & $\ba_5$ & $\ba_7$ & $\ba_{11}$ & $\ba_{13}$ & 
$\ba_{17}$ & $\ba_{19}$ \\
\hline
$353$ & $[1,1,1,-2,16]$ & $2$ &  & $2$ & $1$ & $1$ & $2$ & $2$ & $0$ \\
\hline
$1059$ & $[1,1,1,-66,-270]$ & $2$ &  & $2$ & $1$ & $1$ & $2$ & $2$ & $0$ \\
\hline
$3177$ & $[1,-1,0,-594,6691]$ & $1$ & & $1$ & $1$ & $2$ & $2$ & $1$ & $0$ \\
 & $[1,-1,0,-63,-176]$ & $1$ & & $1$ & $1$ & $2$ & $2$ & $1$ & $0$ \\
\hline
$9531$ & $[0,0,1,3,4]$ & $1$ & & $1$ & $1$ & $2$ & $2$ & $1$ & $0$ \\
 & $[0,0,1,-87891,-10029164]$ & $1$ & & $1$ & $1$ & $2$ & $2$ & $1$ & $0$ \\
 & $[0,0,1,27,-115]$ & $2$ &  & $2$ & $1$ & $1$ & $2$ & $2$ & $0$ \\
 & $[0,0,1,-791019,270787421]$ & $2$ &  & $2$ & $1$ & $1$ & $2$ & $2$ & $0$ \\
\hline
$28593$ & $[1,-1,1,-2162,-38150]$  & $2$ &  & $2$ & $1$ & $1$ & $2$ & $2$ & $0$ \\
 & $[1,-1,0,-240,1493]$  & $1$ & & $1$ & $1$ & $2$ & $2$ & $1$ & $0$ \\
\hline
$85779$ &  & & & & & & & & \\   
\hline
\end{tabular}
\end{center}
Note there are no elliptic curves of conductor $3^5 \cdot 353$.
We see that all the elliptic
curves of conductor $3^k \cdot 353$ have
$a_2 \neq 0 \mod 3$. Thus we are done.
$\qed$
\end{Proof}

\subsection{$p = 5$}

\begin{theorem} Let $E$ be the elliptic curve
$[1,0,1,-80,-275]$ of conductor $7 \cdot 67$.
Let $\rhobar: \Gal(\Qbar/\Q) \rightarrow \GL_2(\F_5)$ be the
representation induced from the $5$-torsion of $E$.
Then $\rhobar$ is surjective, $N(\rhobar) = 67$,
and  $\rhobar$ does not arise from
any elliptic curve of conductor $5^k \cdot 67$.
\end{theorem}

Since $E$ is semistable, the mod $5$ representation is
either reducible or surjective, and a quick check eliminates
the first possibility. Since $E$ is semistable at $7$, we
may check that $E[5]$ is unramified at $7$ by considering
the $7$-adic valuation of the minimal discriminant.
The minimal discriminant is $\Delta_E = -7^{5} \cdot 67$,
and since $5 | 5$ we conclude that $N(\rhobar) = 67$.
Specifically it arises from a modular form of weight
$2$ and level $67$, (in this case
coming from part of the $5$ torsion on a modular
abelian variety $A_f$ of dimension two).

\begin{Proof}
The proof is easier for $p=5$ than for $p=3$, since the largest power
of $5$ dividing the conductor of elliptic curve is two.
Thus we simply enumerate the elliptic curves of conductor
$67$, $67 \cdot 5$, and $67 \cdot 5^2$, and check using
mod $5$ congruences that none of the mod $5$ representations
give rise to $\rhobar$, since $\mathrm{Trace}(\rhobar(\Frob_2)) = 1 \mod 5$.
$\qed$
\end{Proof}

\begin{center}
\begin{tabular}{|c|l|c|}
\hline
$N$ & $E$ & $a_2$ \\
\hline
$67$ & $[0,1,1,-12,-21]$ & $2$ \\
\hline
$335$ & $[0,0,1,-2,2]$ & $0$ \\
\hline
$1675$ & $[0,0,1,-50,281]$ & $0$ \\
 & $[0,-1,1,-13,23]$ & $0$ \\
     & $[0,-1,1,-308,-1982]$  & $3$\\
   & $[0,1,1,-333,2244]$ & $0$\\
\hline
\end{tabular}
\end{center}

\section{Large $p$}

We conclude with a few remarks about $p \ge 7$.
Let $\rhobar: \Gal(\Qbar/\Q) \rightarrow \GL_2(\F_{p})$ be
a modular Galois representation of level $\Gamma_0(N)$ and weight $2$.
Elliptic curves with Galois representations corresponding to
$\rhobar$ are classified by non-cuspidal rational points on  the twisted
modular curve $X(p)(\rhobar)$.
If $ p \ge 7$, then $X(p)$ has
genus $> 1$, and thus there are only finitely many elliptic
curves $E$ which give rise to $\rhobar$, and typically
one would not necessarily expect there to be any.
If $f$ is an eigenform of weight $2$ and level $\Gamma_0(N)$ with
coefficients \emph{not} in $\Q$, then one would expect the mod $\ell$
reductions for $p \ge 7$ also to typically not arise from an elliptic curve.

\begin{theorem} \label{theorem:generic}
 Let $f \in S_2(\Gamma_0(N))^{new}$,  let the coefficients
of $f$ generate the ring $\Or_f$, and assume that $\Or_f \ne \Z$.
 Then for all but finitely
many primes $\p$ of $\Or_f$, the representation
$$\rhobar_{\p}: \Gal(\Qbar/\Q) \rightarrow \GL_2(\Or/\p) 
\hookrightarrow \GL_2(\Fbar_p).$$
does not come from an elliptic curve. 
\end{theorem}

\begin{Proof} Let $f \in S_2(\Gamma_0(N))^{new}$, and let $\Or_f$
be the ring generated by coefficients of $f$. 

\begin{lemma} If $\Or_f \ne \Z$,
then there exists a prime $\ell$ with $(\ell,N) = 1$  and such
that $a_{\ell}(f) \notin \Q$.
\end{lemma}

\begin{Proof} Let $\p$ be a prime in $\Or_K$ of
residue characteristic $p$ such that
$\Or_K/\p \ne  \F_p$. Suppose moreover that $(p,2N)=1$.
Then the associated Galois representation
$$\rho_{\p}: \Gal(\Qbar/\Q) \rightarrow \GL_2(\Or_K/\p)$$
has image that does not land within $\F_p$. This is because the
Hecke eigenvalues $a_q$ for $q | N$ are automatically $0$ or $1$
since $f$ is a newform, and the eigenvalue $a_p$ for $p > k =2$
is determined from the mod $p$ representation~\cite{Edix}.
Now the fact that the trace of $a_{\ell} \mod p$ is the
trace of Frobenius for $(\ell,Np) = 1$ and the fact that
Frobenius elements are dense guarantees an infinite number
of such primes $\ell$. $\qed$
\end{Proof}

Note the lemma and theorem are not true for oldforms.
Take $N = 33$. Then the old form of level $11$ has coefficients
in $\Q(x)/(x^2 + x + 3)$, yet $a_{\ell} \in \Z$ for all $\ell \ne 3$.
Moreover, the mod $p$ representation coming from these old forms
is exactly the mod $p$ representation coming from the elliptic curve
$X_0(11)$.

\medskip

Fix such an $\ell$ as in the lemma above.
 Now suppose that the mod $p$ representation
attached to $f$ comes from an elliptic curve $E$. Assume
that $p \ge 5$. If
$E$ has additive reduction at $\ell$ then since $p \ge 5$,
the field $\Q(E[p])$ is ramified at $\ell$ with ramification index divisible
by $2$ or $3$. This forces the Serre conductor of $\rhobar$ to
be divisible by $\ell^2$ which forces $a_{\ell}$ to be zero,
contradicting our assumption that $\ell \not\in \Q$.
If $E$ has good reduction, then
 $a_{\ell}$ is determined by the mod $\ell$ representation,
and satisfies the Hasse bound $-2 \sqrt{\ell} < a_{\ell}(E) <
2 \sqrt{\ell}$. Moreover $a_{\ell} \equiv a_{\ell}(E) \mod p$.
If $E$ has multiplicative reduction at $\ell$, then either $\ell | N$
in which case $a_{\ell} = \pm 1$ (which is impossible
if $a_{\ell} \not\in \Q$), or one
can ``raise the level'' in the sense of Ribet~\cite{Ribet}.
This is possible only if
$a^2_{\ell} \equiv (1 + \ell)^2 \mod p$. In particular, if
$$A(\ell) = (a^2_{\ell} - (1 + \ell)^2) \prod_{|i| < 2 \sqrt{\ell}}
 (a_{\ell} - i),$$
then  $A(\ell) \equiv 0 \mod p$. Yet $A(\ell)$ is independent of $p$, and
since $a_{\ell}$ is not in $\Q$, $A(\ell)$ is non-zero, and
thus there are only finitely many such $p$. Note in any example
we may explicitly rule out all but finitely many primes $\p$. $\qed$
\end{Proof}

\begin{theorem} \label{theorem:one}  Let $p \ge 11$. Then there exists a modular
semistable Galois representation $\rhobar: \Gal(\Qbar/\Q)
\rightarrow \GL_2(\F_p)$ of weight $2$ and level $\Gamma_0(N)$
that does not arise from any elliptic curve.
\end{theorem}

\begin{Proof} We apply explicitly the proof of theorem~\ref{theorem:generic}
to three particular forms: the form  $f \in S_2(\Gamma_0(23))$ with
$\Or_f = \Z[(1 + \sqrt{5})/2]$, the form $f \in S_2(\Gamma_0(39))$ with
$\Or_f = \Z[\sqrt{2}]$ and the form $f \in S_2(\Gamma_0(590))$ with
$\Or_f = \Z[\sqrt{10}]$. 
For example, for $f \in S_2(\Gamma_0(23))$ we may take
$\ell = 2$, since $a_2
= (\sqrt{5} - 1)/2$.  Then $N_{K/\Q}(A(2))$ is divisible by only
the primes $5$ and $11$. Thus the associated representations for
$p > 5$ and $p \ne 11$  do
not come from elliptic curves. Moreover, the representation
has image inside $\GL_2(\F_p)$ whenever the prime $p$ splits in $\Or_f$,
or equivalently whenever
$(5/p) = 1$. Similar calculations for
the other forms show that if $p >  11$, we are done whenever
$(2/p) = 1$ or $(10/p) = 1$.  Yet since
$$\left(\frac{2}{p}\right)\left(\frac{5}{p}\right) = 
\left(\frac{10}{p}\right),$$
not every quadratic residue symbol can
equal $-1$,
and thus we have found a $\GL_2(\F_p)$ representation
that does not arise from an elliptic curve
for all $p > 11$.
For  $p=11$ one can use the same idea, except with
a different form, for example
$f \in S_2(\Gamma_0(62,\F_{11}))$
with $a_3 \equiv 6 \mod 11$. $\qed$
\end{Proof}

\medskip

This proof does not apply to $p = 7$, since all
possibilities for $\mathrm{Trace}(\rhobar(\Frob_2)) \mod 7$
arise from elliptic curves.
If $a_2 \equiv \pm 3 \mod 7$, then $E$ must necessarily
be semistable at $2$ but have a
$7$-torsion module that is unramified at $7$ (and
so by Tate's theory necessarily have $a_2 \equiv \pm (1 + 2)
\mod 7$). 
For example, when $N = 55$,
there is a form $f$ with coefficients in $\Z[\sqrt{2}]$ and
$a_2 = 1 + \sqrt{2}$. Composing this with the reduction map
to $\F_7$ that sends $a_2$ to $4 \mod 7$ we obtain a candidate
$\rhobar$. Raising the level, we see this representation occurs
from the mod $7$ reduction of a newform of level $2 \cdot 55 = 110$.
Indeed, we find that the $7$ torsion on the  elliptic curve
$E:=[1,0,1,-89,316]$  gives rise to $\rhobar$.
(An easy way to check this is to compute that
$\Delta_E = - 2^7 \cdot 5 \cdot  11^3$, and so the associated mod $7$ representation
is unramified at $2$ and so comes from a mod $7$ representation
of level $55$.) Nevertheless, we prove the following:

\begin{theorem} \label{theorem:two}
There exist irreducible representations $\rhobar :\Gal(\Qbar/\Q)
\rightarrow \GL_2(\F_7)$ with cyclotomic determinant 
that do not arise from elliptic curves over $\Q$.
\end{theorem}

\begin{Proof}
We construct such representations directly. Let $F = \Q(\sqrt{-7})$, and
let $p$ be an inert prime in $\Or_F$ such that $p \equiv -1 \mod 16$ (for example,
$p = 31$). Inside the
ray class field over $F$ of conductor $(p)$ there exists a Galois (over $\Q$)
extension $K$ of degree $8$ over $F$ such that $\Gal(K/\Q)$ is dihedral,
and $K/F$ is totally ramified at $p$.
Let $L = \Q(\zeta_7)$, and let $H = K.L$. Then $\Gal(H/F)$ is cyclic of
degree $24$,
and $\Gal(F/\Q)$ acts on $\Gal(H/F) \simeq \Z/24 \Z$ as multiplication by
$7$ (fixing the subgroup of order $3$, and as inversion on the subgroup of
order $8$). Thus there is a map
$$\rhobar: \Gal(H/\Q) \hookrightarrow \GL_2(\F_7)$$
with image of index exactly two inside the normalizer of the non-split Cartan subgroup
(recall the normalizing element acts as conjugation on $\F^{\times}_{49} \simeq \Z/48 \Z$,
and thus as multiplication by $7$).
 A suitable renormalizing of the non-split Cartan ensures that $\rhobar$ has
cyclotomic determinant. Another realization of $\rhobar$ is the $\omega^2$-twist
of the dihedral representation $\eta: \Gal(K/\Q) \rightarrow  \GL_2(\F_7)$
(in this optic we also observe that $\rhobar$ is modular),
where $\omega$ is the cyclotomic character. The
determinant of $\eta$ is the quadratic character of conductor $7$, which is
$\omega^3 \mod 7$. Thus $\omega^2 \otimes  \eta$ has determinant
$\omega^7 = \omega \mod 7$.
Let us now prove that $\rhobar$ does not come from an elliptic curve. Assume
that $\rhobar$ arises from the $7$-division points of $E/\Q$. All elliptic
curves over $\Q$ acquire semistable reduction after an extension of degree at most
$6$. Moreover, for an elliptic curve with semistable reduction at a prime
$p \ne 7$, the action of inertia at $p$ on the $7$-torsion is either
trivial (in the case of good reduction, by N\'{e}ron--Ogg--Shafarevich) or
factors through a cyclic $7$-group (as can be seen from Tate's parameterization).
We see that  $\# \rhobar|_{I_p} = 8$ is incompatible
with either possibility. $\qed$
\end{Proof}

\medskip
If $(N,p)=1$ and $p \ge 5$, the curve $X_0(p^2 N)$ acquires semistable reduction
over an extension of degree $(p^2 - 1)/2$ (see~\cite{Edix2}). Presumably
many of the $\F_7$-representations of this level have significant inertia
at $p$, and thus  do not arise from elliptic
curves for the  reasons  above. It would be interesting, however, to
find an example of an irreducible
 representation $\rhobar$ such that $X(7)(\rhobar)$ has
points over \emph{every} local completion of $\Q$ but no rational points.
We note that Dieulefait~\cite{Die}  has
also proved theorems~\ref{theorem:one} and \ref{theorem:two},
the first using similar arguments with a rational
form of weight $4$, and the second by finding a representation
whose local representation at $2$ does not come from
an elliptic curve.

\noindent \it Email address\rm:\tt \  fcale@math.harvard.edu

\end{document}